\documentclass[12pt]{article}
\usepackage{latexsym}
\usepackage{amsmath}
\usepackage{amssymb}
\usepackage{amssymb, amscd}
\usepackage{ifthen}
\usepackage{amsthm}
\input epsf
\usepackage[T1]{fontenc}
\usepackage[latin1]{inputenc}

\newtheorem{dfn}{Definition}[section]

\newtheorem{lem}[dfn]{Lemma}

\newtheorem{prop}[dfn]{Proposition}
\newtheorem{cor}[dfn]{Corollary}

\newtheorem{ques}[dfn]{Question}

\def\tupl#1{\langle #1\rangle}
\def\b{\mathbf}
\def\s{\mathsf}
\def\d{\boldsymbol\delta}
\def\l#1{\vskip3pt\noindent\bf Lemma \the\secNum.#1\sl.}
\newcount\secNum\relax\secNum=0
\def\Sec{\global\advance\secNum by 1\relax\the\secNum. }

\def\R{{\Bbb R}}
\def\H{{\Bbb H}}
\def\Z{{\Bbb Z}}

\def\S{{\Bbb S}}

\def\N{{\Bbb N}}

\def\P{{{\cal P}}}

\def\al{\alpha}

\def\ga{\gamma}
\def\Th{\Theta}
\def\g{\gamma}
\def\Ga{{\cal G}}
\def\G{{\Gamma}}

\def\d{\delta}

\def\ve{\varepsilon}

\def\De{\Delta}

\def\th{\theta}

\def\La{\Lambda}

\def\d{\delta}

\def\ti{\widetilde}

\def\nc{\s{Nc}}

\def\a{\mathsf{(a)}}
\def\b{\mathsf{(b)}}

\def\odg{\overline{\d}_g}
\def\od1{\overline{\d}_1}

\def\isom{{\rm Isom}{\H^n}}
\def\iso4{{\rm Isom}{\H^4}}

\def\h3{{\Bbb H}^3}
\def\h4{{\H^4}}
\def\sp{\S^{n-1}_\infty}
\def\sp3{{\S}^{3}_\infty}

\def\bx{$\hfill\square$}

\def\d{\delta}
\def\bu{{\bf u}}

\def\of{F}

\def\BF{\partial_fG}

\def\OG{\overline {G}_f}

\def\1GF{\G_1\cup \partial_f\G_1}

\def\sm0{{\rm Small}_0^{\ve}}
\def\sm1{{\rm Small}_1^{\ve}}

\def\act{\curvearrowright}

\def\DG{\partial_fG}

\def\isom{{\rm Isom}{\H^n}}
\def\iso3{{\rm Isom}{\H^3}}

\def\hcp{{\s H_c p}}
\def\hc{{\s H_c}}
\def\lf1{L_{f,1}}
\def\lfv{L_{f,v}}
\def\bsn{{\mathsf S}^nT}
\def\bs2{{\mathsf S}^2T}
\def\bden{{\bf \Delta}^nT}
\def\bde2{{\bf \Delta}^2T}

\begin{document}

\title{Quasi-isometric maps and Floyd boundaries of relatively hyperbolic groups.}

\author{Victor Gerasimov and Leonid Potyagailo\footnote{{\it 2000
Mathematics Subject Classification.}
  Primary 20F65, 20F67 - Secondary 30F40, 57M07, 22D05
 \hfill\break Key words: Floyd boundary, convergence actions, quasi-isometric maps. }}
\date{\today}
\maketitle \noindent

\bigskip
\bigskip

\begin{abstract}

 We describe the kernel of the canonical map from the
Floyd boundary of a relatively hyperbolic group to its Bowditch
boundary.

Using the Floyd  completion we further prove that the property of
relative hyperbolicity is invariant under quasi-isometric maps. If
a finitely generated group $H$ admits a quasi-isometric map
$\varphi$ into a relatively hyperbolic group $G$  then $H$ is
itself relatively hyperbolic with respect to  a system of
subgroups whose image under $\varphi$ is situated within a
uniformly bounded distance from the parabolic subgroups of $G$.

We then generalize the latter result to the case when $\varphi$ is
an $\al$-isometric map for any polynomial distortion function
$\al.$

As an application of our method we provide in the Appendix  a new
short proof of a basic theorem of Bowditch characterizing
hyperbolicity.

\end{abstract}
\bigskip
\bigskip



\section{Introduction.}

 \noindent We study the actions of discrete groups  by
 homeomorphisms of compact Hausdorff  spaces such
 that
\medskip

 \noindent $\mathsf{(a)}$ the induced action on the space of triples
of distinct points is properly discontinuous, and

\noindent  $\mathsf{(b)}$ the induced action on the space of pairs
of distinct points is cocompact.

\medskip

 Let $T$ be a compact Hausdorff topological space
(compactum). Denote by  $\Th^n T$  the space of subsets of
cardinality $n$  of $T$ endowed with the natural product topology.

Recall that an action of a discrete group $G$  by homeomorphisms
on $T$ is said to have  {\it convergence} property if it satisfies
the condition $\mathsf{(a)}$  \cite{Bo2}, \cite{Tu2}. We also say
in this case that the action of $G$ on $T$ is {\it
$3$-discontinuous}.

An action of $G$ on $T$ is called {\it cocompact on pairs} or {\it
$2$-cocompact} if $\Th^2 T/G$ is compact.

It is shown in \cite{Ge1} that an action  with the properties
$\mathsf{(a)}$, $\mathsf{(b)}$ is {\it geometrically finite} that
is every limit point is either conical or bounded parabolic. From
the other hand it follows from [Tu3, Theorem 1.C] that every
 minimal geometrically finite action on a metrizable
compactum has  the property $\mathsf{(b)}$.

An action of $G$ on $T$ is said to be {\it parabolic} if it has a
unique fixed point. The existence of a non-parabolic geometrically
finite action of a finitely generated group  $G$ on some compactum
$T$ is equivalent to the relative hyperbolicity of $G$ with
respect to some collection of proper subgroups \cite{Bo1},
\cite{Ya}.

  W.~Floyd \cite{F} introduced the notion of
 a boundary    of a finitely generated group as follows.
The word metric   of $G$ is rescaled by a "conformal factor"
$f:\N\to \R_{\geq 0}$ satisfying the  conditions (3-4) below (the
function $f$ we further call  {\it Floyd function}). The Cauchy
completion of the resulting metric space is called {\it Floyd
completion} and is denoted by $\overline G_f$. The
 {\it Floyd boundary}  is the subspace $\partial_fG={\overline G_f}\setminus G$.
The action of $G$ on itself by left multiplication extends to a
convergence action  on  $\overline G_f$ \cite{Ka}.

It is shown in \cite{F} that for any geometrically finite discrete
subgroup $G<\iso3$ of the isometry group of the hyperbolic space
$\Bbb H^3,$ and for  a quadratic scaling  function $f,$ there
exists a continuous $G$-equivariant map $F$ from $\DG$ to the
limit set $T=\La(G)$. The preimage $F^{-1}(p)$ of a point $p$ is
not a single point if and only if $p$ is a parabolic point of rank
1 in which case it is a pair of points \cite{F}. P.~Tukia
generalized Floyd's Theorem to geometrically finite discrete
subgroups of $\isom$ \cite{Tu1}.

 If an action of  a finitely generated group $G$
  on a
 compactum $T$ has the properties    $\a$ and $\b$
  then by \cite{Ge2} there exists a continuous equivariant map $F_0:\partial_{f_0} G\to  T$ for the exponential
    function $f_0(n)=\lambda^n$ where $\lambda\in ]0,1[$. An easy
    argument then shows (see Corollary  \ref{polynmap} below) that
    a continuous $G$-equivariant map $F: \partial_fG\to T$ exists
    for any Floyd function $f$ such that $f\geq f_0$.  So the function $f_0$
    is a critical function starting from which the whole theory
    works. The map $F$ obtained in this way   we
    call   {\it Floyd map}.

The group $G$ is a discrete  open subset of $\overline G_f$. For a
subset $H$ of $G$ denote by $\partial H$ its topological boundary
in the space $\overline G_f$. Let $\s{Stab}_{G} p$ denote  the
stabilizer in $G$
 of a point $p\in T$. Since the action of $G$  has the convergence property
 the topological boundary   of $\s{Stab}_{G} p$ coincides with the limit set
$\Lambda(\s{Stab}_{G} p)$ for the action of $\s{Stab}_{G} p$ on
${\overline G}_f.$
\par Our first
result describes the kernel of the map $F$.

\bigskip

  \noindent {\bf Theorem A.} {\it Let
  $G$ be    a finitely generated  group acting   on a
  compactum $T$   $3$-discontinuously and $2$-cocompactly.
  For every Floyd function $f$ such that $f\geq f_0$ let $F:\DG\to T$
  be the Floyd map. Then
   $$F^{-1}(p)=\partial(\s{Stab}_{G} p)\hfill\eqno(1)$$
   for any parabolic point $p\in T$. Furthermore, $F(a)=F(b)=p$   if and only if either $a=b$ or $p$ is
parabolic.}

\bx

\bigskip

Note that  the subgroup inclusion does not imply an
  embedding of Floyd boundaries so we have:

\begin{ques}

Let  a finitely generated group $G$ act $3$-discontinuously and
$2$-cocompactly on a compactum $T$. Let $F : \DG\to T$ be a
continuous $G$-equivariant map. For a parabolic point $p\in T$ is
it true that
$$\partial (\s{Stab}_{G}p)=\partial_{f_1}\s{Stab}_{G}p$$

\noindent  for  some scaling function $f_1 ?$
\end{ques}

In particular, an affirmative answer to this question would give a
complete generalization of the Floyd theorem for relatively
hyperbolic groups.

 Our next result describes the quasi-isometric
 maps into relatively hyperbolic groups:

\bigskip

\noindent {\bf Theorem B.} {\it Let a finitely generated group $G$
act  $3$-discontinuously and $2$-cocompactly on a compactum $T$.
Let $\varphi :H\to G$ be a quasi-isometric map of a finitely
generated group $H$ into $G$. Then there exist   a compactum $S$,
a $3$-discontinuous $2$-cocompact action of $H$ on $S,$ and a
continuous map $\varphi_*:S\to T$ such that  for every
$H$-parabolic point $ p{\in}S$  the point $\varphi_* p$ is
$G$-parabolic, and $\varphi(\s{Stab}_H p)$ is contained in a
uniformly bounded neighborhood  of a right coset of
$\s{Stab}_G(\varphi_* p)$ in $G.$\bx

}
\bigskip

\noindent Using known facts about the relative hyperbolicity
Theorem B can be reformulated as follows:

\begin{cor}
\label{rel}

Let $G$ be a finitely generated relatively hyperbolic group with
respect to a collection  of subgroups $P_j\ (j=1,...,n)$. Let $H$
be a finitely generated group and  let $\varphi : H\to G$ be a
quasi-isometric  map. Then $H$ is relatively hyperbolic with
respect to a collection  $Q_i$ such that $\varphi$ maps each $Q_i$
into a uniform neighborhood of a left coset  of some $P_j$ (the
case $Q_i=H$ is allowed).\bx
\end{cor}

\bigskip

  Note that several particular cases of
Corollary \ref{rel} were already known: if either  the
quasi-isometric map $\varphi : H\to G$ is a monomorphism [Hr,
Corollary 1.6]; or   $\varphi$ admits a quasi-isometric inverse
map [Dr, Theorem 1.2]; or  $G$ is relatively hyperbolic with
respect to subgroups which themselves are not relatively
hyperbolic with respect to proper subgroups [BDM, Theorem 4.8].
Theorem B provides a complete control of parabolic subgroups of
the initial group $H$ under the map $\varphi$ (compare with [BDM,
page 546]). The proof is independent and does not rely on the
above results.

After our preprint has appeared on the ArXiv \cite{GePo1}, we were
informed by the authors of \cite{BDM} that even though Corollary
\ref{rel} is not explicitly stated in their paper, the proof of it
can be obtained using several results of [Dr] and [BDM].

\bigskip

We call a scalar non-decreasing function $\al:\N\to \R_{>0}$ {\it
distortion function} if $\al_n=\al(n)\geq n\ (\forall n\in\N)$.

\begin{dfn}
\label{aldist} A map $\varphi :X\to Y$ between the metric spaces
$(X,d_x)$ and $(Y, d_Y)$ is called {\it $\alpha$-isometric} if
$$d_Y(\varphi(x), \varphi(y))\leq\al(d_x(x,y))\ {\rm and}\
d_X(x,y)\leq \al(d_Y(\varphi(x), \varphi(y)))\hfill\eqno(1)$$

\noindent for a distortion function $\al.$
\end{dfn}

\medskip

\noindent Note that the partial cases of $\al$-isometric maps are:
1) $\al(n)=n$ (isometric); 2) $\al(n)=cn$ (Lipschitz); 3)
$\al(n)=cn+c$ (quasi-isometric).

Our next goal is to generalize Theorem B to the case when the
distortion function of the map $\varphi$ is polynomial.

\bigskip

\noindent {\bf Theorem C.} {\it Let a finitely generated group $G$
act  $3$-discontinuously and $2$-cocompactly on a compactum $T$.
Let $\varphi :H\to G$ be an $\al$-isometric map of a finitely
generated group $H$ into $G$ for a polynomial distortion function
$\al.$ Then all statements of Theorem B are true for the group $H$
and the map $\varphi.$}\bx

\bigskip

 We now  outline the
content of the paper. In Section 2 we provide some preliminaries.
In Section 3 we prove several important lemmas implying  the proof
of Theorem A. The proof of Theorem B is given in Section 4. In
Section 5 we first prove that every $\alpha$-geodesic   which is
far from the origin   has a small Floyd length (Generalized
Karlsson Lemma). Then following the lines  of the proof of Theorem
B we prove Theorem C.

Even though Theorem B is a partial case of Theorem C we prefer to
keep both statements in the paper by two reasons. Historically we
first obtained the proof of Theorem B appeared in the earlier
version of our paper \cite{GePo2}. Furthermore since the
quasi-isometric maps are themselves very important,  Theorem B has
its independent interest. So we first give a complete proof of
Theorem B and then
   develop all necessary tools giving the proof of Theorem C.

  As an application of our methods we give in
the Appendix a short proof of the fact that the existence of a
$3$-discontinuous and 3-cocompact action on a compactum without
isolated points implies that the group is word-hyperbolic
\cite{Bo3}.

\medskip
{\bf Acknowledgements.} During the work on this paper both authors
were partially supported by the ANR grant ${\rm BLAN}~07-2183619.$
We are  grateful to the Max-Planck Institute f\"ur Mathematik in
Bonn, where a part of the work was done. We also thank the
Brasilian-French cooperation grant having  supported our work.

The authors are very much thankful to the referee for her/his
valuable suggestions for the improvement of the paper

\section{Preliminaries.}

\subsection{Convergence actions.} By {\it compactum} we   mean a
compact Hausdorff space. Let $\bsn$ denote the orbit space of the
action of the permutation group on $n$ symbols on the  product
$\displaystyle T^n=\underbrace{T{\times}\dots{\times}T}_{\mbox
{$n$ times}}$. The elements of $\bsn$ are {\it the generalized
non-ordered $n$-tuples}. Let $\Th^nT$ denote the open subset of
$\bsn$ consisting of the non-ordered $n$-tuples whose components
are distinct. Put $\bden=\bsn\setminus\Th^nT$. So
  $\bde2$ is   the image of the diagonal of $T^2$.

\bigskip

\noindent {\bf Convention.} {\rm If the opposite is not stated all
group actions on compacta are assumed   to   have the convergence
property. We will also assume that $\vert T\vert > 2,$ and so
$\Th^3T\not=\emptyset.$}

\bigskip

Recall a few common definitions  (see e.g. \cite{Bo2}, \cite{GM},
\cite{Fr}, \cite{Tu2}).  The {\it discontinuity domain}
$\Omega(G)$ is the set of the points of $T$ where $G$ acts
properly discontinuously. The set $\Lambda (G)=T\setminus
\Omega(G)$ is the {\it limit set} and  the points of $\Lambda(G)$
are called {\it limit points}. An action of $G$ on $T$ is called
{\it minimal} if $\La(G)=T.$

It is known that $\vert\La(G)\vert\in\{0,1,2, \frak c\}$
\cite{Tu2}.  An  action    is called {\it elementary} if its limit
set   is finite.

 A point $p\in T$ is called {\it parabolic} if $\vert\Lambda({\rm
 Stab}_Gp)\vert=1$.

A limit point $x\in \La(G)$ is called {\it conical} if there
exists an infinite sequence of distinct elements $g_n\in G$ and
distinct points $a, b\in T$ such that
$$\forall y\in T\setminus\{x\}\ :\ g_n(y)\to a\in T\ \wedge g_n(x)\to
b.$$ Denote by $\nc T$ the set of non-conical points of $T$.

 A
parabolic   point $p\in \La(G)$ is called {\it bounded parabolic}
if the quotient space $(\La(G)\setminus \{p\})/\s{Stab}_G p$ is
compact.

 An action  of $G$ on $T$  is called  {\it
geometrically finite} if every non-conical limit point  is
bounded parabolic.

A subset $N$ of the set $M$ acted upon by $G$ is called {\it
$G$-finite} if it intersects finitely many $G$-orbits.

\begin{lem}
\label{parfin} [Ge1, Main Theorem] If the action of $G$ on $T$ is
$3$-discontinuous and $2$-cocompact then

\begin{itemize}

\item[a.] The set $\nc T$ is $G$-finite.

\medskip

\item[b.] For every $p\in\nc T$ the quotient $(T\setminus p)/{\rm
Stab_pG}$ is compact. \bx
\end{itemize}
\end{lem}

\medskip

It follows from Lemma \ref{parfin} that for a $2$-cocompact action
of $G$ on $T$  a non-conical point $p\in \nc T$ is isolated in $T$
   if and only if its stabilizer $\s{Stab}_pG$ is finite. Hence, a non-conical
   point with infinite stabilizer   is
   bounded parabolic.

\subsection{Quasigeodesics and Floyd completions of graphs.}

Recall that a {\it ($c$-)quasi-isometric map}   $\varphi : X\to Y$
between two metric spaces $X$ and $Y$ is an $\al$-isometric map
for the affine distortion function $\al(n)=cn+c\ (c > 1):$

$${1\over c}d_X(x, y) - c \leq d_Y(\varphi(x),\varphi(y))\leq cd_X(x, y)+ c,\hfill\eqno(2)$$

\medskip
\noindent where $d_X,\ d_Y$ are the metrics of $X$ and $Y$

\bigskip

\noindent {\bf Remarks.} 1) {\rm A quasi-isometric map can in
general be   multivalued. This more general case can be easily
reduced to the case of a one-valued map.

2) Sometimes   other terms for quasi-isometric maps are used: {\it
large-scale Lipschitz maps} \cite{Gr1} or {\it quasi-isometric
embeddings} \cite{BH}} (note that $\varphi$ is not necessarily an
injective map).

\bigskip

 A \it path \rm is a distance-nonincreasing
map $\gamma:I\to \G$ from a nonempty convex subset $I$ of $\Bbb
Z$. The \it length \rm of a path  $\gamma$ is the diameter of $I$
in $\Bbb Z$. A \it subpath \rm is a path which is a restriction of
$\gamma$.

A path $\gamma:I\to \G$ is called $c$\it-quasigeodesic \rm if it
is   a $c$-quasi-isometric  map.   In the case when $\ga$ is an
isometry, it is called  \it geodesic\rm.

 A  ($c$-quasi-)geodesic path $\ga : I\to \G$ defined on a
half-infinite subset $I$ of $\Bbb Z$ is called

 \noindent {\it ($c$-quasi-)geodesic ray}; a (quasi-)geodesic path
defined on the whole $\Bbb Z$ is called \it ($c$-quasi-)geodesic
line\rm.

Let
  $d(,)$ denote  the canonical shortest path distance function
on $\G$. We denote by $\s N_DM$ the $D$-neighborhood of a set
$M{\subset}\G$ with respect to $d.$

We now briefly recall the construction of the Floyd  completion of
a
  graph $\G$ due to W.~Floyd   \cite{F}.
  Let  $\G$ be a
locally finite connected graph endowed with a basepoint
$v\in\G^{0}$. Let $f:\Z_{>0}\to \R_{>0}$ be a function satisfying
the following conditions:

\begin{itemize}

\item[] $$\exists\ {K > 0}\ \forall n\in\N\ :\ 1 \leq {f(n)\over
f(n+1)} \leq K\hfill\eqno(3)$$

\medskip

\item[] $$\displaystyle \sum_{n\in \N} f(n)
<\infty.\hfill\eqno(4)$$

\end{itemize}

\bigskip

For  convenience  we extend the function $f$ to $\Z_{\geq 0}$ by
putting $f(0):=f(1).$ We further call the function $f$ satisfying
(3-4) {\it Floyd (scaling) function}.

Define the {\it Floyd length} of an edge joining vertices $x$ and
$y$ as $f(n)$ where $n=d(v, \{x,y\}).$ Then the Floyd  length
$L_{f,v}$ (or simply $L_v$ for a fixed Floyd scaling function $f$)
  of a path is the sum of the Floyd lengths of its
  edges. The Floyd distance $\d_v=\d_{f,v}(a,b)$ is
 the shortest path distance:

 $$\displaystyle \d_{v}(a,b)=\inf_{\alpha} \lfv(\al),\hfill\eqno(5) $$

\noindent where the infimum is taken over all paths between  $a$
and $b.$

It follows from  (3) that every two metrics $\d_{v_1}$ and
$\d_{v_2}$ are bilipschitz equivalent with a Lipschitz constant
depending on $d(v_1, v_2).$   The Cauchy completion
${{\overline\G}_f}$
  of the metric space $(\G^0, \d_v)$   is
  called {\it Floyd completion}. It
  is compact and does not depend on
the choice of the basepoint $v.$  Denote by $\partial_f\G$ the set
${\overline\G}_f\setminus \G$ and call it {\it Floyd boundary}.
The distance $\d_v$ extends   naturally   to ${\overline\G}_f$.

The following Lemma   shows that the Floyd length of a far
quasigeodesic  is small.

\begin{lem}
\label{agraph} {\rm (Karlsson Lemma)} For  every $\varepsilon
> 0$ and  every $c >0,$  there exists  a finite set $D\subset \G $ such that $\d_v$-length
of every  $c$-quasigeodesic  $\ga\subset \G$   that does not meet
$D$  is less than $\ve$.\bx

\end{lem}

{\bf Remark.}   {\rm    A.~Karlsson \cite{Ka} proved it for
geodesics in the Cayley graphs of finitely generated groups. The
proof of \cite{Ka} does not use the group action and is still
valid for quasigeodesics.}

\medskip

 Let $S$ be a set of paths of the form
$\ga:[0,n[\to \G$
  of unbounded length starting at the same point $a=\al(0).$
    Every path $\ga\in S$  can be considered as an element of the product
$\displaystyle \prod_{i\in I} \s N_{ i}(a)$. Since $\G$ is a
locally finite graph the space $\displaystyle\prod_{i\in I} \s N_{
i}(a)$ is compact in the Tikhonov topology. So every infinite
sequence of paths in $S$ contains a
 subsequence converging to a path of the form
$\delta:[0,+\infty[\to \G$ all of whose initial segments are
initial segments of paths  in $S$. In particular if $S$ is a set
of $c$-quasigeodesics (or more generally $\alpha$-geodesics, see
Section 5) then every infinite subset of $S$ admits a subsequence
converging to a $c$-quasigeodesic (respectively $\alpha$-geodesic)
of infinite length.
 Note that the infinite
limit path exists in a more general context when $\delta(0)$
belongs to a fixed finite set.

\begin{dfn}
\label{geodray} {\rm  For a $c$-quasigeodesic ray
$r:[0,\infty[\to\G$ we say that $r$ {\it converges to a point} in
$\partial_f\G$ if the sequence $(r(n))_n$ is a Cauchy sequence for
the $\d_f$-metric. We also say in this case that $r$ {\it joins}
the points $r(0)$ and $\displaystyle x=\lim_{n\to\infty}r(n)\in
\partial_f\G$.}
\end{dfn}

\medskip

\begin{prop}
\label{geodline}

 Let $\G$ be a locally finite connected graph. Then

 \begin{itemize}

  \item[\sf a.] For each $c>0$ every $c$-quasigeodesic ray in $\G$ converges
   to a point in $\partial_f\G$.

\medskip

  \item[\sf b.] For every  $p{\in}\partial_f\G$ and every $a\in \G$
there exists a geodesic ray joining $a$ and $p$.

\medskip

  \item[\sf c.] Every two distinct points in $\partial_f\G$ can be joined by a
geodesic line.

\end{itemize}
\end{prop}

\proof {\sf a.} Let $r:[0, \infty[\to\G$ be a $c$-quasigeodesic
ray. Put $x_n=r(n)$ and $r_n=r([n,\infty[)$. For any vertex $v\in
\G^{0}$ we have $d(v, r(n))\to\infty.$ It follows  from Karlsson
Lemma that $L_{f, v}(r_n)\to 0.$

\medskip

{\sf b.} Let $B_f(p,R)$ denote the ball in the Floyd metric   at
$p\in\partial_f\G$ of radius $R.$ For $n{\geqslant}1,$ choose
$a_n\in  B_f(p,{1\over n})$ and join $a$ with $a_n$  by a geodesic
segment $\gamma_n$. Let   $\gamma$ be the limit path for the
family $S = \{\gamma_n:n{>}0\}$. By  ({\sf a})  $\ga$ converges to
a point $q{\in}\partial_f\G$. If $p\ne q$ set $3\delta = \d_1(p,q)
>0$ where $1\in \G$ is a fixed vertex (we use this notation
keeping in mind the case of  Cayley graphs).   Let $n$ be an
integer for which
 $L_{f, 1}(\gamma\vert_{[n,\infty[})\leq \delta$. For $m\geqslant
n$ we can choose $k$  such that $\gamma_k|_{[0,m]} =
\gamma|_{[0,m]}$ and $\d_1(a_k, p){\leqslant}\delta$. So $L_{f,
1}(\gamma_k\vert_{[m,k]})\geq \delta$. Since the distance $d(1,
a_n)$ is unbounded, by Karlsson Lemma the quantity $L_{f,
1}(\gamma_k\vert_{[m,k]})$ should tend to zero. This contradiction
shows that $p=q.$

\medskip

{\sf c)} Let $p,q{\in}\partial_f\G$ and $p{\ne}q$. By  $({\sf b})$
there exist geodesic rays $\alpha,\beta:[0,\infty[\to \G$ such
that $\al(0)=\beta(0)=a$ and $\al(\infty)=p,\ \beta(\infty)=q.$
Let  $3\delta = \d_v(p,q) $. By Karlsson Lemma every geodesic
segment joining a point $\alpha(n)$ in $ B_f(p,\delta)$ with a
point $\beta(n)$ in $B_f(q,\delta)$ intersects a finite set
$B=B(a, R)\subset\G$. So there exists an infinite sequence of
geodesic segments $\gamma_n$ passing through a point $b\in B$
whose endpoints converge to the pair $\{p,q\}$. A limit path for
such sequence is a geodesic line in question.  \bx

\bigskip

 Let $\G_i\ (i=1,2)$   be  locally finite connected graphs.
  The following Lemma gives a sufficient condition
to extend a quasi-isometric map between the Floyd completions of
the graphs.

\begin{lem}
\label{homfl} Let  $\varphi : \G_1\to\G_2$ be a
$c$-quasi-isometric map for some $c\in \N.$ Suppose that there
exists a constant $D>0$ such that

$${f_2(n)\over f_1(c n)} < D\ (n\in \N)\hfill\eqno(6),$$

\noindent  Then  the map $\varphi$ extends   to a Lipschitz map
between the Floyd completions ${\overline \G}_{1, f_1}
\to{\overline \G}_{2, f_2}$.

 \end{lem}

\proof To simplify the notations put ${\overline \G}_{i,
f_i}={\overline \G}_i\ (i=1,2).$    We denote by $E$ and $K$ the
constants  from (3)  corresponding to  $f_1$ and $f_2$
respectively such that $f_1(n)\leq E\cdot  f_1(n+1)$ and
$f_2(n)\leq K\cdot f_2(n+1)\ (n\in \N).$ By omitting the  indices
we use the notations  $d$ and $\d$ for  the graph distances
 and the Floyd distances   with respect to the chosen
basepoints  denoted by $1$ in  both graphs $\G_i\ (i=1,2).$

  We first prove that
$\varphi$ is a Lipschitz map with respect to the Floyd metrics,
i.e.

$$\forall\ x,y\in \G_1\ : \ \d(x, y) \geq \ve\d(\varphi (x), \varphi (y))\hfill\eqno(7)$$

\noindent for some  $\ve>0$.

It suffices to prove the statement for the case when $d(x,y) = 1$.

By  (2) we have $d(\varphi (x), \varphi (y)) < cd(x,y) +c=2 c.$
Let $\ga : [0,n]\to \G_2$
 be a geodesic realizing the distance $d(\varphi(x), \varphi(y)),$
and let $a_i=\ga(i)\ (i=0,...,n)$  be its vertices where
$a_0=\varphi (x), a_n=\varphi (y).$ We have
$\d(x,y)=f_1(d(1,\{x,y\})).$ Assume that $d(1,\{x,y\})=d(1,x).$
Then $$d(1, a_i)\geq d(\varphi (1), \varphi (x)) - d(\varphi (x),
a_i)-d(1, \varphi(1))\geq  d(\varphi (1), \varphi (x)) -2c-
d(1,\varphi(1)) = d(\varphi (1), \varphi (x)) - n_0,$$ where
$n_0=2c+ d(1,\varphi(1)).$

 Let $r_0$ be a number such that $r_0 >  c(n_0+c)$ and $B(1,r_0)$ denotes
 the ball
 centered at $1$ of the radius $r_0$. Then for
 $x\in\G_1\setminus B(1, r_0)$  we have $d(\varphi(1), \varphi(x)) > n_0$.
 Then using the
monotonicity of $f_2$ and   (3) for $x,y\in \G_1\setminus
B(1,r_0)$ we obtain:
$$\displaystyle \d(\varphi (x), \varphi (y))\leq\sum_{i=0}^{n-1} f_2(d(1,
\{a_i,a_{i+1}\}))\leq \sum_{i=0}^{n-1} f_2(d(\varphi(1), \varphi
(x))-n_0)\leq  2c K^{n_0}f_2(d(\varphi(1),\varphi(x))).$$

\noindent Since $d(1,x) > c^2$ the last term can be estimated
using (2), (3) and (6) :

$$\displaystyle f_2(d(\varphi(1),\varphi(x)))\leq Df_1(c d(\varphi(1),\varphi(x)))\leq
Df_1(c\cdot d(1,x)/c - c^2)\leq D\cdot E^{c^2}\cdot f_1(d(1,x)),$$

 \noindent Summing all up we
conclude

$$\displaystyle\d(\varphi (x), \varphi (y))\leq
 2cDK^{n_0}E^{c^2}\cdot f_1(d(1,x)).$$

So (7) is true for the constant $\displaystyle
\ve=(2cDK^{n_0}E^{c^2})^{-1}$  outside of the ball $B(1, r_0)$. By
decreasing the constant $\ve$ if necessary, we obtain the
inequality (7) everywhere on $\G_1.$

The map $\varphi:(\G_1,\d)\to (\G_2,\d)$  being Lipschitz extends
to a Lipschitz map ${\overline \G}_1 \to{\overline \G}_2$.
 \bx

\bigskip

\noindent {\bf Remarks.} {\rm If, for a function $f$, the value
${f(n)\over f(2n)}$ is bounded from above  (e.g. for any
polynomial function) then one can take the same scaling function
$f_1=f_2=f$ for both graphs $\G_1$ and $\G_2$ independently of
$c$.

\medskip

If  the scaling function for the graph $\G_2$ is $f_2(n)=\al^n\
(\al\in ]0,1[)$ then to satisfy (6) we  can take $f_1(n)=\beta^n$
 as the scaling function for  the group
$H$ where $\beta = \al^{1/c}$}. \bx

\bigskip

Let $G$ be a finitely generated group and let $S$ be  a   finite
generating set for $G$.  Denote by $d$ the word metric
corresponding to $S.$ Let $\OG$ denote the Floyd completion
$G\sqcup\partial_fG$ of the Cayley graph of $G$ with respect to
$S$  corresponding to a
 function $f$ satisfying   (3-4). Condition
    (3)   implies the equicontinuity for the action of every element $g\in G$
     by left multiplication  on $G$ (with an equicontinuity constant depending on $g$) \cite{Ka}.
     So the  action of each element $g\in G$ extends by homeomorphism to
     $\OG$. Therefore the  whole group $G$ acts on $\OG$ by homeomorphisms.  The Floyd
metric $\d_g$ is the $g$-shift of $\d_1$ (where $1$ is the
identity element of $G$):

$$\d_g(x, y) =\d_1(g^{-1}x, g^{-1}y),\ \ x,y\in\OG, g\in G.$$

\bigskip

\noindent On the space $\OG$ we   also consider  the following
  shortcut   pseudometrics.

\begin{dfn}\label{psmetric} \rm Let $\omega$ be a closed $G$-invariant equivalence relation
on $\OG.$  The {\it shortcut pseudometric} $\overline{\d}_g$ is
the maximal element in the set of symmetric functions $\varrho:
\OG{\times}\OG\to\Bbb R_{\geqslant0}$ that vanish on $\omega$ and
satisfy the triangle inequality, and  the inequality
$\varrho{\leqslant}\d_g$.
\end{dfn}

For $p,q{\in}\OG$ the value $\odg(p,q)$ is the infimum of the
finite sums $\displaystyle\sum_{i=1}^n\d_g(p_i, q_i)$ such that
$p{=} p_1$, $q{=}q_n$ and $\tupl{ q_i, p_{i+1}}{\in}\omega\
(i{=}1,\dots,n{-}1)$ [BBI, pp 77]. Obviously,    $\odg$  is the
$g$-shift of the pseudometric $\od1.$  The pseudometrics
$\overline{\d}_{g_1}$ and $\ \overline{\d}_{g_2}$ are bilipschitz
equivalent with the same constant as for $\d_{g_1}$ and
$\d_{g_2}$.

The pseudometric $\odg$   induces a  pseudometric on the quotient
space $\OG/\omega.$ We denote this induced pseudometric by the
same symbol  $\odg$.

\medskip

 The following result will be often used further.

\begin{lem}
\label{Flmap}\cite{Ge2} Let $G$ be a finitely generated group
acting $3$-discontinuously and $2$-cocompactly   on a compactum
$T\ (\vert T\vert >  2).$    Then there exists
 a constant   $\lambda\in]0, 1[$ and a continuous $G$-equivariant map
 $F_0:\partial_{f_0} G\to T$
  where $f_0(n)=\lambda^n$.
\end{lem}

\medskip

 \noindent {\it Summary of the proof of Lemma \ref{Flmap}.} For the reader's convenience
  we provide few explanations about the proof of the existence of the Floyd map. We
 refer to the preprint \cite{Ge2} for more details. It
 consists of two parts. In the first part it is proven a general
 statement [Ge2, 5.2]  that if $G$ is a locally
 compact group acting on a compactum  $T$ and on a locally compact Hausdorff space
 $L$ such that  the action on $T$
 is 3-proper and the action on $L$  is proper and
 cocompact, then there exists a unique topology on the compactified
 space $T\cup L$. Furthermore the induced topologies on $T$ and $L$
 coincide with the original ones and the union of two actions
 gives rise to a $3$-proper action $G\act (T\cup L).$

To prove the Lemma  it is enough however to consider a partial
case  of the latter statement when a finitely generated group $G$
acts 3-discontinuously on the Hausdorff space $T$. In this case
one can directly show that there exists a locally finite,
connected $G$-finite graph $\Ga$ such that the action of $G$ on
the topological space $M=T\cup\Ga$ is 3-discontinuous [GePo2,
Proposition 3.14]. The graph $\Ga$  is a graph   of entourages of
$T$ (i.e. neighborhoods of the diagonal $\Delta^2 T$ of $T^2$) and
the above finiteness properties of $\Ga$ are proved in \cite{Ge1},
\cite{GePo2}.

Since $G$ acts 2-cocompactly on $T$ it is easy to prove that the
action $G\act M$   also admits a compact fundamental set (see e.g.
Lemma \ref{GAT} below). Then by [Ge1, Proposition E] there exists
a generating entourage $\bu$ on $M$, i.e. the orbit $G\bu$
generates the set of entourages $\mathsf{Ent} M$ of $M$ as a
filter. Then it follows from [Ge2, 6.2] that the system
$\{\bu_n=\cap(S^n\bu)\ \vert \ n\in \N\}$ is a {\it Frink system}
of entourages, where $S$ is a finite symmetric generating set of
$G$ and $S^n$ is a collection of words of $G$ of length at most
$n.$ By a classical result from the general topology   the above
Frink system gives rise to a Frink metric ${ \Delta}_{\bu, S}$ on
$M$ which is the maximal among all the metrics $\rho$ on $M$
satisfying $\rho\vert_{\bu_n}\leq 2^{-n}.$ Furthermore if
${\Delta}_{\bu, S}(x,y)\leq 2^{-n}$ then $(x,y)\in \bu_{n-1}$ [Ke,
Lemma 6.2] (or [Ge2, Proposition 6.1.1]). In particular the Frink
metric generates the topology of $M$  and $\{\bu_n\ \vert\
n\in\N\}$ is a countable base for the uniformity of this topology.
Let  $v\in \Ga^0$ be  a basepoint.  It is shown in [Ge2, 6.3] that
there exists  an exponential  Floyd function $f(n)=2^{-n}/\rho$
such that for the Floyd metric $\delta_{f,v}$, and
 for every edge $e\in \Ga^1$ one has
$$\Delta_{\bu, S} (e) \leq c\delta_{f,v}(e),\hfill\eqno(*)$$

\noindent where       $\rho$ and $c$ are uniform constants not
depending on $e$. So (*) implies that the inclusion map $\Ga^0\to
M$ is uniformly continuous. Thus since $M$ is complete, the
identity map on $\Ga$ extends continuously and equivariantly to
the  map $F:{\overline \Ga}_{f}\to M$ satisfying

$$\forall x, y\in{\overline
\Ga}_f\ \ :\ \  \Delta_{\bu, S} (F(x), F(y)) \leq
c\delta_{f,v}(x,y),\hfill\eqno(**)$$

\noindent  Since the graph $\Ga$  is $G$-finite and connected
  there exists a $G$-equivariant
$C$-quasi-isometry $\psi$ between the Cayley graph of $G$ (with
respect to some finite generating set) and the graph $\Ga$. By
Lemma \ref{homfl} it extends to a $G$-equivariant map $\psi :
{\overline G}_{f_0}\to {\overline \Ga}_{f}$ where $f_0=f^{1/C}$
and $f(n)=2^{-n/\rho}$ is the above Floyd function for $\Ga.$ So
$f_0(n)=\lambda^n$ where  $\lambda=2^{-1/C\rho}$. The map
$F_0:=F\circ \psi$ satisfies the statement of the Lemma. \bx

\medskip

\noindent {\bf Remark.} The same proof gives that
 for any connected locally compact
graph $\G$ admitting a cocompact action of the group $G$  there
exist an exponential function  $f(n)=\mu^n\ (\mu\in ]0,1[)$ and a
$G$-equivariant continuous  map $F:
\partial_f\G\to T$.
\bx

\medskip

The Floyd function $f_0$ given by Lemma \ref{Flmap} is a critical
function satisfying the following.

\begin{cor}\label{polynmap}
 For every   function  $f\geq f_0$ satisfying (3-4)
 there exists a continuous $G$-equivariant  map $F:\partial_fG\to T$. In particular
 this is true for any   function
   $f(n)=1/\P(n)$ where $\P(n)$ is a
polynomial  of degree $k>1\ (n\in\N\cup\{0\},\ \P(0)\not=0)$.
\end{cor}

\proof Let $f_0(n)=\lambda^n$ be the function from \ref{Flmap}.
Since $f(n)\leq f_0(n)\ (n>n_0)$ we have $\delta_{f_0}(x,y)\leq
\delta_{f}(x,y)\ (x,y\in G).$ So  the identity map extends
continuously and equivariantly to the map $\chi:\partial_{f}
G\to\partial_{f_0} G$. By Lemma \ref{Flmap} there exists a
continuous $G$-equivariant  map $F_0:\partial_{f_0}G\to T$. Put
$F=F_0\circ \chi.$

Since every function of type $1/\P$ where $\P$ is a polynomial of
degree $k>1$ satisfies the conditions (3-4) and $f_0(n)\leq
1/\P(n)\ (n>n_0)$ it admits such a map too. \bx

\medskip

\noindent {\bf Definition.}  The map $F$ given by Corollary
\ref{polynmap}   we   call {\it Floyd map.}

\section{The orbit compactification  space $\ti T$ and its convex subsets.}

In this Section we fix a $3$-discontinuous and $2$-cocompact
action by homeomorphisms of a finitely generated group $G$ on a
compactum $T$ containing at least 3 distinct points.

\subsection{The space $\ti T$.} For  a fixed Floyd scaling function
$f$ let $F:\partial_fG\to T$  be the Floyd map obtained in Lemma
\ref{Flmap}. We extend $F$ over $\OG=G\sqcup\partial_fG$ to the
disjoint union $\ti{\ti T}= T\sqcup G$  by the identity map ${\rm
id} : G\to G$. We keep the notation $F$ for this extension. The
maps $T\overset{\s{id}}\to\widetilde {\ti
T}\overset{F}\leftarrow\OG$ determine on $\widetilde {\ti T}$ the
pushout topology: a set $S{\subset}\widetilde {\ti T}$ is open if
and only if $S{\cap}T$ is open in $T$ and $F^{-1}S$ is open in
$\OG$. The space $\widetilde {\ti T}$ being the union of two
compact spaces $T$ and $F(\OG)$ is a compactum.

Since for every point $x\in\partial_fG$ is the limit of a sequence
$g_n\in  G$ and the map $F$ is equivariant we obtain that
$\Lambda(G)= F(\partial_f\G)$ for the limit set $\La(G)$ of the
action $G\act T.$ By Lemma \ref{parfin}   the set $T\setminus
\Lambda(G)$ is $G$-finite. Denote by $\ti T$ the subspace
$\La(G)\sqcup G$ of $\ti{\ti T}.$

\medskip

\noindent {\bf Remark.} We need to introduce $\ti{\ti T}$ before
  $\ti T$ in order to include the exceptional case of 2-ended groups.
In this case $\La(G)$ consists of 2 points and we need at least
one more point to apply Lemma \ref{Flmap}.

\medskip

\begin{prop}
\label{eqvmap} Let $G$ act on  compacta $X$ and $Y$ and let $\psi
:X\to Y$ be a  $G$-equivariant continuous surjective map. If the
action of $G$ on $X$ is $3$-discontinuous, then the action of $G$
on $Y$ is $3$-discontinuous.
\end{prop}

\proof  The
  map $\psi$
induces a proper $G$-equivariant continuous surjective map
${\mathsf S}^3X\to {\mathsf S}^3Y$.  Let $K$ and $L$ be compact
subsets of $\Th^3Y.$ Since $Y$ is Hausdorff the preimage of every
compact  in $\Th^3 Y\subset {\mathsf S}^3Y$ is compact in $\Th^3
X$. Thus   $ K_1=\psi^{-1}(K)$ and $L_1=\psi^{-1}(L)$ are compact
subsets of $\Th^3 X$. The action on $X$ is discontinuous so the
set $\{g\in G\ \vert\ g  K_1\cap L_1\not=\emptyset\}$ is finite.
By the equivariance of   $\psi$    the set $\{g\in G\ \vert\
gK\cap
 L\not=\emptyset\}$ is finite too. \bx

\medskip

 \begin{lem}
 \label{GAT}
 The induced action $G$ on $\ti T$ is  $3$-discontinuous and
$2$-cocompact.
\end{lem}

\proof  By \cite{Ka} the group $G$ acts $3$-discontinuously on
$\OG=G\sqcup \DG$. The Floyd map $F:\OG\to\ti T$ is
$G$-equivariant and continuous. By  Proposition \ref{eqvmap}  the
action on $\ti T$ is $3$-discontinuous.

If  $K$ is a compact fundamental set for the action of $G$ on
$\Th^2(T)$ then $K_1=K\cup \{1\}\times(\ti T\setminus\{1\})$ is a
compact fundamental set for the action of $G$ on $\Th^2 \ti T.$\bx

\bigskip

 Let $\omega$ be the kernel of the Floyd map
$F:\OG\to \ti T,$ i.e. $\displaystyle (x,y)\in\omega$ if and only
if $F(x)=F(y)$. It determines the shortcut pseudometric
$\overline{\d}_g\ (g\in G)$ on $\OG$ (see Subsection 2.2). The map
$F$  transfers it to  a pseudo-metric on $\ti T$ also denoted by
$\overline{\d}_g$:

$$\forall x,y\in\BF\ \ :\ \overline{\d}_g(F(x), F(y))=\overline{\d}_g(x,y).\hfill\eqno(8)$$

\noindent Since $\overline{\d}_g$ is the maximal pseudometric on
$\OG$ satisfying $\overline{\delta}_g\leq\delta_g$ for the Floyd
metric  $\delta_g$, by the property (**) of Section 2 the
transferred pseudometric $\overline{\d}_g\ (g\in G)$ becomes a
real metric on $\widetilde T$, i.e.
$$\forall p,q\in \ti T\ :\ \overline{\d}_g({p,q}) = 0\Longrightarrow\ p = q.\hfill\eqno(8')$$

\noindent Any metric $\overline{\d}_g$  determines the topology of
$\widetilde T$.

\begin{lem}
\label{finpar} Let $H$ be the stabilizer of a parabolic point $p$.
Every $H$-invariant set $M{\subset}G$ closed in $\widetilde
T{\setminus}\{p\}$ is $H$-finite.\end{lem}

\proof By Lemma \ref{GAT} the action of $G$ on  $\ti T$ is
$3$-discontinuous and $2$-cocompact, so by Lemma \ref{parfin} the
space
  $(\widetilde T{\setminus} \{p\})/H$ is compact.
Since $G\subset\ti T$ is an orbit of isolated points, the closed
subset $M/H$ of $(\widetilde T{\setminus} \{p\})/H$   consists of
isolated points. So the set $M/H$ is finite.\bx

\bigskip

\subsection{Horocycles and horospheres}

  By Proposition \ref{geodline}.a  every $c$-quasigeodesic
ray $\gamma:\Bbb N\to G$ converges to a point $ p\in \partial_fG$.
We call the point $\displaystyle p=\lim_{n\to\infty}\ga(n)$   \it
target \rm of $\gamma$ and denote it by $\gamma(\infty)$. The path
$F\circ\gamma$ converges to the point $F(p)\in\ti T$ which we also
call \it target\rm. In other words a $c$-quasigeodesic ray extends
to a continuous map from $\overline{\Bbb N}=\Bbb
N{\cup}\{{+\infty}\}$ to $\widetilde T$. The target is necessarily
a limit point for the action $G\act \ti T.$ So any  bi-infinite
$c$-quasigeodesic $\gamma:\Bbb \Z\to G$ extends to a continuous
map of ${\overline \Z}=\Z{\cup}\{{-}\infty,{+}\infty\}$ with
$\gamma({\pm}\infty){\subset}T$. The set  $\gamma(\pm\infty)$ is
either a pair of limit points or a single limit point.

\medskip

Every  quasigeodesic, either finite or infinite, is defined on a
  closed  subset $J{\subset}\overline Z$ with $\partial J$
being a pair or a single point.

\begin{dfn}
\label{horocycle} {\rm A bi-infinite $c$-quasigeodesic
$\gamma:\Z\to G$ is called {\it  $c$-horocycle} at $p\in T$ if
 $\ga(+\infty)=\gamma(-\infty)=p.$}

\end{dfn}

\medskip

\begin{dfn}
\label{hull}

{\rm
 The $c$\it-hull \rm$\s H_c M$  of  a set $M{\subset}\widetilde T$    is the union of $M$ and
  all $c$-quasigeodesics (finite or infinite)  having the endpoints in
$M$:
 $$\displaystyle{\s H_cM}=M{\cup}\{\gamma(I)\ \vert\ \gamma:I\to G\
{\rm is\ a\ c-quasigeodesic},\ I\subset \Z,\ {\rm and}\
\gamma(\partial I) \subset M\}.$$

\noindent The $c$-hull $\s H_cp$ of  a   single-point set $\{
p\}\in T$ is called $c$\it-horosphere at $ p$.}

\end{dfn}

\medskip

\noindent By $\overline M$ we denote the closure of $M$ in
$\widetilde T$.

\medskip

 \noindent{\bf Main Lemma.} {\it $T \cap \overline
M=T \cap \overline{\s H_c M}$ for every    $M\subset\ti T$ and
  $c>0$.}

\proof  Suppose by contradiction that there exist $M$ and  $c$
such that  $T \cap \overline{\s H_c M}\setminus\overline
M\ne\emptyset$. Let $a\in T \cap \overline{\s H_c
M}\setminus\overline M.$ By Karlsson Lemma there exists a number
$r>0$ such that the $\d_1$-length of every $c$-quasigeodesic
outside of the ball $N_r(1)\subset G$ is less than $\ve={\overline
\d}_1(\overline M, a)/2>0$ where $1\in G$ is a basepoint.

\bigskip

 \centerline{\epsfxsize=3.5in\epsfbox{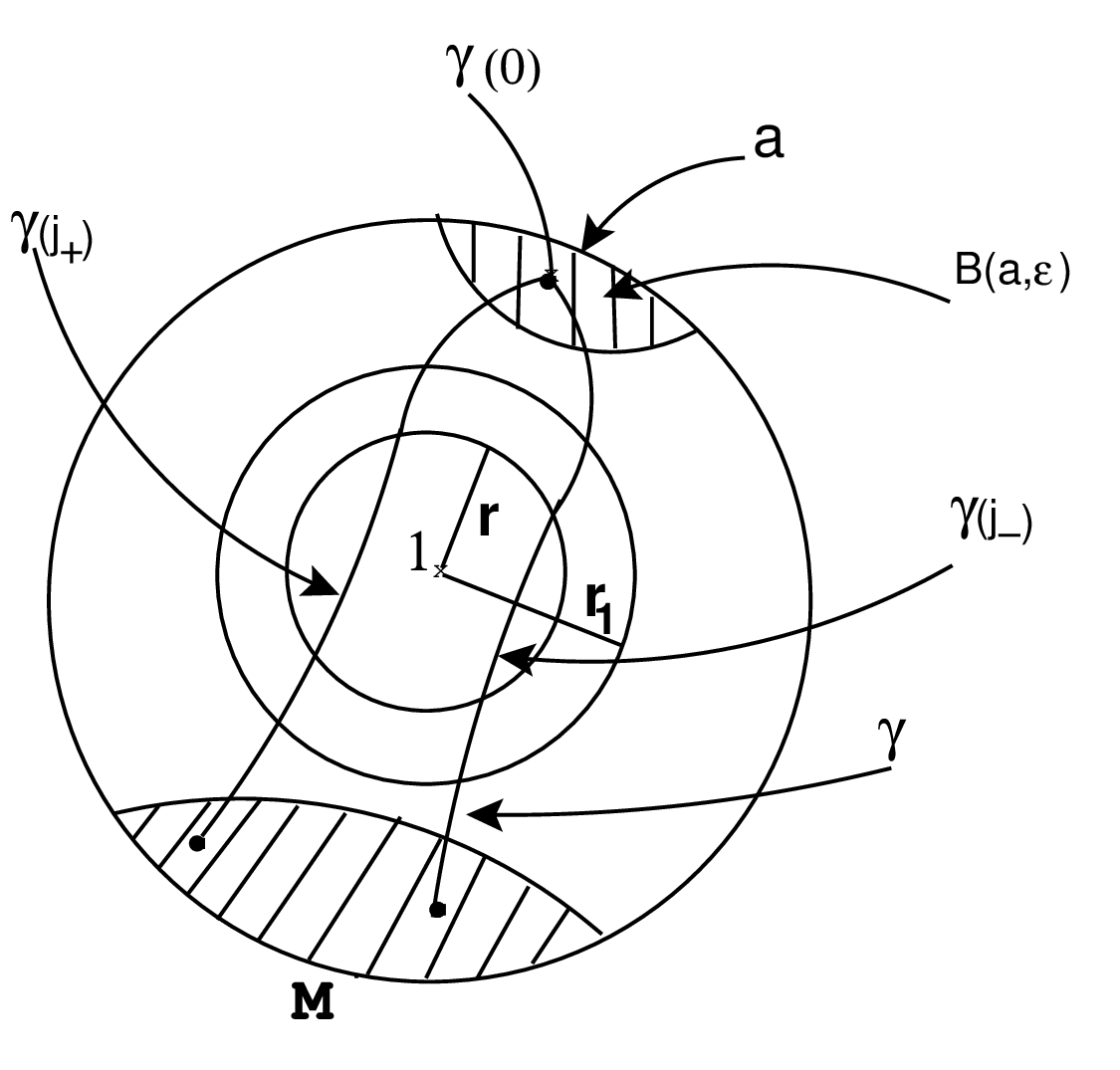}}

By the assumption there exists a $c$-quasigeodesic $\ga: I=[i_-,
i_+]\to G$ such that $0\in I,\  \ga(\partial I)\subset M,$ and
$\ga(0)$ is arbitrarily close  to $a.$ So we can assume that
 ${\overline\d}_1(a,\gamma(0)) < \ve $ and   $ \gamma(0)\not \in N_{r_1}(1)$
where $ r_1=r+cr+{c\over 2}$. Let $\gamma_+=\gamma|_{I\cap\Bbb
N}$, $\gamma_-=\gamma|_{I\cap(-\Bbb N)}$. For their Floyd lengths
$L_1$ (see 2.2) we have
$$L_1(\gamma_\pm) \geq {\overline\d}_1(\ga(i_{\pm}),\ga(0))\geq
{\overline\d}_1(\ga(i_{\pm}), a) -{\overline\d}_1(a,\ga(0))\geq
{\overline\d}_1(\overline M, a){-}\overline{\d}_1(a,\ga(0))\geq
2\ve - \ve = \ve.$$

\noindent So there exists a subseqment $J = [j_-,j_+]$   of $I$
such that $0{\in}J$ and $d(1,\gamma(j_\pm)){\leqslant}r$.

We obtain $d(\ga(j_-), \ga(j_+))=\s{diam} (\gamma(\partial
J)){\leqslant}2r$ and $\s{length}(\gamma|_J){\leqslant}c(2r+1)$.
Thus $d(1,\gamma(0))\leq d(1, \gamma(\partial J))+d(\gamma(0),
\gamma(\partial J))\leq r_1=r+ {c\over 2}(2r+1).$ So $\gamma(0)\in
N_{r_1}(1).$  A contradiction.\bx

\bigskip

\noindent We finish the subsection by obtaining several
corollaries of the Main Lemma.

\medskip

\begin{lem}
\label{nohorcon}
 There is no $c$-horocycle  at a  conical
point.\end{lem}

\proof  Let $p\in  T$ be a conical point. There exist distinct
points $a, q\in T$ and a sequence $(g_n)\subset G$ such that
$g_n(p)\to q$ and $g_n(x)\to a$ for all $x\in  T\setminus\{p\}.$
 Suppose  by contradiction  that  $\gamma$ is
an $c$-horocycle at $ p$.

Let $Q$ be a closed neighborhood of $q$ such that $ a\notin Q$. We
can assume that $ q_n=g_n (p)\in  Q$ for all $n$. So
$g_n(\gamma(0))\in \s H_c Q$. By Lemma \ref{GAT} the action
$G\act\ti T$ has   the convergence property. So $g_n(\gamma(0))\to
a$. It follows that $a\in (\overline{\s H_c Q}\cap T)\setminus
(Q\cap T)$ contradicting to the Main Lemma.\bx

\medskip

  Until the end of this Subsection we fix a parabolic fixed point
   $p$, and denote   by $H$ the stabilizer of $p$ in $G.$

\begin{lem}
\label{hulfin} For every $c$ the set $(G \cap \s H_c p)/H$ is
finite.
\end{lem}

\proof The set $\s H_c p$ \rm is $H$-invariant. By the Main Lemma
it is closed in $\ti T{\setminus} p$ and $p$ is its unique limit
point. Thus the set $G \cap \s H_c p$ is a closed $H$-invariant
subset of $\widetilde T{\setminus} p$. The result   follows from
Lemma \ref{finpar}.\bx

\medskip

\begin{lem}
\label{finfin}
 The closure in $\ti T$ of any $H$-finite subset $M$ of $G$ is
$M{\cup}\{ p\}$.\end{lem} \proof It suffices to consider the case
when $M$ is an $H$-orbit. As $d(M,\s H_c p)$ is bounded, the Floyd
distance $\d_1(m,\s H_c p)$ tends to zero while $m\in M$ tends to
$T$. \bx

\bigskip

\begin{cor}
\label{parqconv} There exists a constant $C_0$ such that the
stabilizer of every parabolic point is  $C_0$-quasiconvex.
\end{cor}

\proof Let $\hcp$ be a $c$-horosphere at $p$. By Lemmas
\ref{hulfin} and \ref{finfin} the set $\hcp\cap G$ is $H$-finite.
So $d(\hcp, H)\leq D$ and $H\subset M$ where $M$ is a
$D$-neighborhood   of $\hcp\cap G.$

We have $\overline M=M\cup \{p\}.$ By the Main Lemma
$\overline{\hc M}=\hc M\cup \{p\}.$ So $\hc M$ is closed in $\ti
T\setminus \{p\}$ and by Lemma \ref{finpar} it is also $H$-finite.

Let $\ga:I\to G$ be a geodesic segment  with   endpoints in
$H\subset M.$ Then $\ga(I)$ and $M$ are both  subsets of the
$H$-finite set $\hc M.$ Hence for any $a\in\gamma(I)$ there exist
$h_i\in H\ (i=1,2)$ and $b\in M$ such that $d(h_1(a), h_2(b))\leq
{\rm const}$. Since $M$ is $H$-invariant we have
$h_1^{-1}h_2(b)\in M$. Then $d(a, M)\leq {\rm const}$ and so $d(a,
H)\leq {\rm const}.$

By Lemma \ref{parfin} the set of parabolic points is $G$-finite so
there exists a uniform constant $C_0$ such that every stabilizer
of a parabolic point is $C_0$-quasiconvex. \bx

 \medskip

\subsection{The kernel of the Floyd map.}

  \noindent {\bf Theorem A.} {\it Let
  $G$ be    a finitely generated  group acting   on a
  compactum $T$   $3$-discontinuously and $2$-cocompactly.
  For every Floyd function $f$ such that $f\geq f_0$ let $F:\DG\to T$
  be the Floyd map. Then
   $$F^{-1}(p)=\partial(\s{Stab}_{G} p)\hfill\eqno(1)$$
   for any parabolic point $p\in T$. Furthermore, $F(a)=F(b)=p$   if and only if either $a=b$ or $p$ is
parabolic.}

\proof  Extending the map $F$ by the identity map  as in
Subsection 3.1 we can suppose that we have continuous and
$G$-equivariant map $F:\OG\to \ti T.$ Denote $H=\s{Stab}_{G} p.$
Let $x\in F^{-1}(p)$. We will show that $x\in\partial H.$ Let
$y\in\partial H.$ If $y=x$ then there is nothing to prove. If not,
then by Proposition \ref{geodline}.c there exists a bi-infinite
geodesic $\gamma$ joining $x$ and $y.$ It is a horocycle in
$\widetilde T$, so $\gamma(\Bbb Z)\subset\hcp.$ By Lemma
\ref{hulfin}, $\gamma(\Bbb Z)$ is contained in
$Hg_1{\cup}\dots{\cup}Hg_l$. By Lemma \ref{finfin} the boundary of
each $H$-coset is $\{ p\}$. So $\displaystyle x=\lim_{n\to\infty}
h_ng_i$ where $i\in\{1,...,l\}$ and $ h_n\in H.$ Since $d(h_ng_i,
h_n)=d(g_i, 1)$ is constant for all $n,$ we have $\d_1(x, h_n)\to
0$ and $x\in
\partial H\subset\partial_fG.$

Assume that $ a \ne b$. Then as above join $a$ and $b$  by a
bi-infinite geodesic $\gamma$.  Then $ \ga$ is an horocycle in
$\ti T$, and by Lemma \ref{nohorcon} the point $p=F(a)=F(b)$ is
parabolic. \bx

\bigskip

\noindent{\bf Remark.} By Theorem A the preimage of a conical
point under an equivariant map is a single point. It is true in a
more general setting without assuming that the action of $G$ on
$T$ is $2$-cocompact [Ge2, Proposition 3.5.2].

\begin{cor}
\label{fltu} In the notation of Theorem A the set
$\partial(\s{Stab}_{G} p)$ is a quotient of the Floyd boundary
$\partial_{f_1}(\s{Stab}_{G} p)$ with respect to some scaling
function $f_1.$

\end{cor}

\proof Since $H=\s{Stab}_{G} p$ is undistorted in $G$ \cite{Ge1},
the inclusion map $H\hookrightarrow G$ is $c$-quasi-isometric for
some integer $c.$ For a given Floyd function $f$ there exists a
Floyd function $f_1$ satisfying conditions (3-4) such that
$f(n)/f_1(cn)$ is bounded from above. By Lemma \ref{homfl} the
inclusion extends to a continuous map ${\overline H}_{f_1}\to
{\overline G}_f.$   It maps $\partial_{f_1}(\s{Stab}_{G} p)$ onto
$\partial(\s{Stab}_{G} p).$\bx

 \bigskip
\noindent{\bf Remarks.}  It follows  from Karlsson Lemma that the
Floyd boundary of a virtually abelian group is either a point or
pair of points. In particular it is true for any discrete
elementary subgroup  of $\isom$.

We have to confess that we do not understand  in the paper
\cite{F} the proof of the fact that the preimage by $F$ of a
parabolic fixed point $p\in\La(G)$ for a geometrically finite
Kleinian group $G<\iso3$ belongs to the boundary of its stabilizer
in ${\overline G}_f$ [F, page 216]. Theorem A and Corollary
\ref{fltu} complete this argument  implying that the preimage of a
fixed point of an elementary parabolic subgroup of a geometrically
finite subgroup of $\isom$  of rank at least $2$ is a point.

\section{Compactification of a quasi-isometric map. Proof of Theorem B.}

Let $G$ and  $H$ be finitely generated groups with  fixed finite
generating sets. We denote the corresponding  word metrics by the
same symbol $d$.

Let $\varphi:H\to G$ be a  $c$-quasi-isometric map. We now choose
two Floyd functions $f_i\ (i=1,2)$ satisfying the condition (6) of
Lemma \ref{homfl}. Let $\overline H_{f_1}$ and $\overline G_{f_2}$
denote the Floyd completions corresponding to   $f_1$ and $f_2$
respectively.

\medskip

 \noindent To simplify   notations  we put
$$X=\partial_{f_2}G,\ Y=\partial_{f_1}H,\  \ti X =X\sqcup G,\ \ti Y=Y\sqcup H.$$

 \noindent  By  Lemma
\ref{homfl} it extends to a uniformly continuous map $\ti Y\to \ti
X$ which we keep denoting by   $\varphi.$ By continuity reason the
inequality (7) of Lemma \ref{homfl} remains valid for this
extension.

\noindent The kernel $\theta_0$ of the composition $\widetilde
Y\overset\varphi\to\widetilde X\overset{F}\to\widetilde T$ is a
closed equivalence relation on $\widetilde Y$.  We have

$$\displaystyle (x,y)\in\theta_0 \ \Longleftrightarrow\
F\varphi(x)=F\varphi(y).\hfill\eqno(9)$$

\noindent The following  equivalence on $\ti Y$ is closed and
$H$-invariant:

$$\theta= \cap \{h\theta_0:h{\in}H\}.\hfill\eqno(10)$$

\noindent  So $\displaystyle (x,y)\in\theta$ if and only if
$\displaystyle (h(x), h(y))\in \theta_0$ for each $h\in H.$

Let $\ti S=\ti Y/\theta$ and $S=Y/\theta$. Denote by $\pi$ the
quotient map $\ti Y\to\ti S$. It is $H$-equivariant. Since $\th$
is closed the space $\ti S$ is a compactum [Bourb, Prop I.10.8].
The open subspace $A=\pi(H)$ of $\ti S$ is an $H$-orbit of
isolated points. The group $H$ acts $3$-discontinuously on $\ti
Y$. Since $\vert \ti S\vert>2$   we have by   Proposition
\ref{eqvmap} that the action of $H$ on $\ti S$ is
$3$-discontinuous.

The map $\varphi$ sends $\theta$  into the kernel of $F$, and
induces a continuous map $\varphi_*:\widetilde S\to\widetilde T$
such that the following diagram is commutative:

$$  \begin{CD}\displaystyle  \ti X @>\of>>  \ti T\\ @A{\varphi}AA
@A{\varphi_{*}}AA\\
 \ti Y @>\pi>>\ti S \\
\end{CD}
$$

\bigskip

\noindent The equivalence $\theta$ determines the shortcut
pseudometrics $\overline{\d}_g$ on   $\ti S$  such that
$$\forall\ {p,q}\in\ti Y, \forall g\in G\ \ :\ \ \overline{\d}_g({p,q})= \overline{\d}_g([p],[q]).\hfill\eqno(11)$$

\begin{prop}
\label{cocomppairs}
 The action of $H$ on $\widetilde S$ is
$2$-cocompact.\end{prop}

\proof By Lemma \ref{GAT}  it suffices to verify $2$-cocompactness
on the space $S$. Let $[{p],[q]}$ be  distinct $\theta$-classes in
$S$. For some $h{\in}H$ and $p,q\in Y$ we have $\tupl{h p,h
q}{\notin}\theta_0$, so $F\varphi(hp){\ne}F\varphi(h q).$ Denote
by $K$   a compact fundamental set for the action of $G$ on
$\Th^2\ti T$, i.e.
$$\Theta^2\ti T = {\cup}\{gK:g{\in}G\}.$$ Let $\delta$ be the infimum of
 the continuous function $\overline{\d}_1|_K$. It  is
strictly positive by  (8'). There exists $g{\in}G$ such that
$\{g^{-1}F\varphi(h p),g^{-1}F\varphi(h q)\}{\in}K$. So
$\overline{\d}_g(F\varphi(h p), F\varphi(h q))\geq\delta$. As
$\d_g\geq \overline{\d}_g$ we also have
 $\d_g(\varphi(h p), \varphi(h q))\geq
 \delta$.
 Let $\gamma$ be a bi-infinite geodesic in $H$ with $\partial\gamma =
\{h p,h q\}$. Since $\varphi$ is $c$-quasi-isometric,
$\varphi(\gamma)$ is contained in a $c$-quasigeodesic in $G$.
 So  by
Karlsson Lemma there exists $r = r(c,\delta)$ such that
$d(g,\varphi(\gamma)){\leqslant}r$. Assume that
$d(g,g_0){\leqslant}r$ for  $g_0 = \varphi(h_0),\ h_0=\ga(0)$.

By the bilipschitz equivalence of the shortcut metrics,
$\displaystyle \overline{\d}_{g_0}(F\varphi(h p), F\varphi(h
q)\geq{\delta\over D(r)}$ for a function $D(r)$ depending only on
$r.$ By  $(8),$
  $\overline{\d}_{g_0}(\varphi(hp),
\varphi(hq)) = \overline{\d}_{g_0}(F\varphi(hp), F\varphi(hq)).$
By Lemma \ref{homfl} we obtain

 $$\overline{\d}_1(h_0^{-1}hp,h_0^{-1}h
q)=\overline{\d}_{h_0}(h p,h q)\geq
\ve\overline{\d}_{g_0}(\varphi(hp),\varphi(hq))\geq
{\ve\delta\over D(r)}=\d_1.$$ Using  (11)  we conclude that the
set
$\{\{s_1,s_2\}{\in}\Theta^2S:\overline{\d}_1(s_1,s_2)\geq\delta_1\}$
is a compact fundamental set for the  action of $H$ on
$\Theta^2S$. \bx

\medskip

\bigskip

\noindent {\bf Theorem B.} {\it Let a finitely generated group $G$
act  $3$-discontinuously and $2$-cocompactly on a compactum $T$.
Let $\varphi :H\to G$ be a quasi-isometric map of a finitely
generated group $H$.

Then there exist   a compactum $S$, a $3$-discontinuous
$2$-cocompact action of $H$ on $S,$ and a continuous map
$\varphi_*:S\to T$ such that  for every $H$-parabolic point $
p{\in}S$   the point $\varphi_* p$ is $G$-parabolic, and
$\varphi(\s{Stab}_H p)$ is contained in a uniformly bounded
neighborhood of a right coset of $\s{Stab}_G(\varphi_* p)$ in $G.$

}
\bigskip

 \proof The space $S$ and the map
$\varphi_*$ are already constructed. We are going to prove that
$\varphi_*$ maps $H$-parabolic points to $G$-parabolic points. Let
$ p$ be a parabolic point for the action of $H$ on $\ti S$ and let
$Q$ be its stabilizer.  By [Ge1, Main Theorem, d] $Q$ is finitely
generated and undistorted in $H$.  So the embedding
$Q\hookrightarrow H$ is  quasi-isometric. Since  $Q$ is infinite
there exists a bi-infinite geodesic $\gamma:\Z\to Q$. Thus
$\gamma$ is a $c$-quasigeodesic in $H$ for some $c$. Moreover
since the set of non-conical points $\nc S$  is $H$-finite the
constant $c$ can be chosen uniformly for all parabolic points $p$.

By Proposition \ref{cocomppairs}   the action of $H$ on $\ti S$ is
$2$-cocompact. By  the Main Lemma   the boundary of $Q$ in
$\widetilde S$ is $\{ p\}$. In particular, $\gamma$ is a
$c$-horocycle.  Since $\varphi$ is quasi-isometric, the path
$\varphi\circ\gamma : \Z\to\ti T$ is a $l$-quasigeodesic for some
uniform constant  $l>0$. The continuity of $\varphi_*$ and the
commutativity of the above diagram  imply  that
$\displaystyle\lim_{n\to\pm\infty}\varphi_*(\ga(n))=\varphi_*(p).$
 Thus $\varphi_*\circ\ga$ is a $l$-horocycle at the point
$\varphi_* p.$ It follows from Lemma \ref{nohorcon} that
$\varphi_* p$ is parabolic for the action of $G$ on $T.$

Every $h{\in}Q$ belongs to a bi-infinite geodesic in $Q$. So by
the above argument we have $\varphi(Q)\subset\s H_l(\varphi_*p).$
 Since $\nc T$ is $G$-finite
there exists a finite set $W$ of $G$-non-equivalent parabolic
points of $G$ such that every parabolic point is $G$-equivalent to
a point in $W.$ So there exists $g\in G$ such that
$g(\varphi_*p)=q\in W.$ We have $g(\s H_l(\varphi_*p))=\s H_lq.$
By Lemma \ref{hulfin} the set $G\cap \s H_lq$ is
$\s{Stab}_Gq$-finite. So there exists a uniform constant $c>0$
such that $\forall q\in W\ \forall y\in \s H_lq\ d(y, \s{Stab}_Gq)
\leq c.$ So $d(g\varphi(Q),  \s{Stab}_Gq) = d(\varphi(Q),
g^{-1}\s{Stab}_Gq)= d(\varphi(Q),
\s{Stab}_G(\varphi_*p)g^{-1})\leq c.$ \bx
\bigskip

\noindent {\it Proof of the Corollary \ref{rel}:} Suppose that $G$
is a finitely generated relatively hyperbolic group with respect
to parabolic subgroups $P_i\ (i=1,...n)$ in the strong sens of
Farb \cite{Fa}. Then by \cite {Bo1} (see also \cite{Hr})   the
group $G$ possesses a geometrically finite $3$-discontinuous
action on a compact metrizable space $X$. It follows from [Tu3,
Theorem 1.C] that  the space $\Th^2 X/G$ is compact. Let $S$ be a
compactum as in Theorem B on which the group $H$ acts
 $3$-discontinuously and $2$-cocompactly. By
\cite{Ge1} this action
 is geometrically finite, the set of parabolic points  is   $H$-finite,
  and their stabilizers are all finitely generated.
   It follows from \cite{Ya} that $H$ is
relatively hyperbolic with respect to the stabilizers $Q_i\
(i=1,...,k)$ of $H$-non-equivalent parabolic points.  By Theorem B
the image of every parabolic subgroup of $H$ by $\varphi$ is
contained in a uniform neighborhood of a right coset of the
corresponding parabolic subgroup of $G$ which is a left coset of
one of $P_j\ (j=1,...,n)$. The Corollary follows. \bx

\section{$\alpha$-isometric rigidity. Proof of Theorem C.}

The goal of this Section is to generalize Theorem B to the case of
$\al$-isometric maps (see (1)).

\medskip

\noindent {\bf Theorem C.} {\it Let a finitely generated group $G$
act  $3$-discontinuously and $2$-cocompactly on a compactum $T$.
Let $\varphi :H\to G$ be an $\al$-isometric map of a finitely
generated group $H$ into $G$ for a polynomial distortion function
$\al.$ Then all statements of Theorem B are true for the group $H$
and the map $\varphi.$}

\bigskip

\noindent The proof will follow the lines of the proof of Theorem
B using  few new facts.

\medskip

\noindent {\bf Definition.} A path $\ga:I\to \Gamma$ in a graph
$\Gamma$ is called $\alpha$-geodesic if the map $\gamma$ is
$\al$-isometric.

\bigskip

 \noindent We start with the following generalization of
Karlsson Lemma \ref{agraph}:

\bigskip

\begin{lem}{(Generalized Karlsson Lemma)}
\label{genkarl} Let $\Gamma$ be a locally finite connected graph
endowed with a basepoint $v\in \Gamma^0$. Let $\al:\N\to \R_{>0}$
and $f:\N\to \R_{>0}$ be respectively distortion and Floyd scaling
function satisfying:

$$\displaystyle\sum_{n\in\N}\al_{2n+1}f_n < +\infty,\hfill\eqno(12)$$

\noindent where $\al_n=\al(n)$ and  $f_n=f(n).$ Then for every
$\varepsilon >0$ there exists $r=r(\ve)$ such that for every
$\al$-distorted path $\ga:I\to \G$  if $d(v,\ga)
>  r$  then for its Floyd length we have $L_{\delta_{f,v}}(\ga) <
\ve.$
\end{lem}

\proof Let $0\in I$
  realizes the distance  $d(v,\g)=d(v, \g(0)) > r.$ For every
$i \in ]r, \max I[$ we denote by $t_i$ the natural number such
that $x_i=\min \{x\in I\ \vert \ \forall\ t \geq x\ :\ d(\g(t),
v)\geq i\}.$ We have $d(v, \g(x_i))=i$. Indeed if not then $d(v,
\g(x_i))> i$ and by the triangle inequality $d(v, \ga(x_i-1)) >
i-1$ and so $d(v, \ga(x_i-1)) \geq i$ which is impossible by the
choice of $x_i.$ So the interval $I$ is subdivided in the
intervals $I_i=[x_i, x_{i+1}[$ such that $\forall t\in I_i\ d(v,
\ga(t))\geq i.$  By the triangle inequality  $d(\g(x_i),
\g(x_{i+1}))\leq 2i+1$. Since $\g$ is $\al$-geodesic we have
$x_{i+1}-x_i\leq \al_{2i+1}.$ For the Floyd length  it yields
$L_{\d_{f,v}}(\g\vert_{[x_i, x_{i+1}]})\leq f(d(v, I_i))\cdot {\rm
length}_d \ga(I_i)=\al_{2i+1} f_i$. Thus
$$\displaystyle L_{f,v}(\g)\leq \sum_{i=r}^k \al_{2i+1}f_i
+\al_{2k}f_k,$$

\noindent where $k<\infty$ only if $I$ is  finite and  $d(v,
\ga(x_k))=d(v, \ga(\max I)).$ By (12) there exists $r=r(\ve)$ such
that $\displaystyle \sum_{i=r}^k \al_{2i+1}f_i <\ve/2$ and
$\al_{2k}f_k\leq \al_{2k+1}f_k<\ve/2$ as $(\al_n)_n$ is a
non-decreasing sequence. The Lemma is proved.\bx

\medskip

\begin{dfn}
\label{adm} We call the couple $(f, \al)$ of a Floyd scaling
function $f$ and a distortion function $\al$ {\it admissible} if
it satisfies (12).
\end{dfn}

\begin{cor}
\label{algeodline} Every $\alpha$-geodesic ray in a locally finite
connected graph $\G$ converges to a point in $\partial_f\G$ if the
pair $(\al, f)$  is admissible.
\end{cor}

\proof It repeats the argument  of Proposition \ref{geodline}
where instead of  Karlsson Lemma one uses Generalized Karlsson
Lemma. \bx

\medskip

 \noindent The following result generalizes Lemma \ref{homfl} to
 the case of
$\al$-isometric maps:

\medskip

\begin{lem}
\label{alhomfl} Let $\G_i\ (i=1,2)$ be locally finite connected
graphs and $\varphi :\G_1\to\G_2$ be an $\al$-isometric map for
some distortion function $\al:\N\to \N.$ Suppose that Floyd
scaling functions $f_i:\N\to \R_{>0}\ (i=1,2)$ satisfy

$${f_2(n)\over f_1(\al_n)}\leq D\ (n\in \N)\hfill\eqno(13)$$

\noindent for some constant $D>0.$ Then $\varphi$ extends to a
Lipschitz map (denoted by the same letter) $\varphi
:\overline{\G}_{1,f_1}\to\overline{\G}_{2,f_2}.$
\end{lem}

\proof We keep the notations of the proof of Lemma \ref{homfl}.
Using the $\alpha$-isometric map $\varphi$ the same argument as in
Lemma \ref{homfl} gives that $\forall x,y\in \G_1$ such that
$d(x,y)=1$ the following is true:

$$\delta_{f_2}(\varphi(x), \varphi(y))\leq
2\al_1K^{n_0}f_2(d(\varphi(1), \varphi(x))\leq
2\al_1K^{n_0}Df_1(\al(d(\varphi(1), \varphi(x)))\\ \leq
2\al_1K^{n_0}D f_1(d(1,x)),$$

\noindent where $n_0=2\al_1+d(1, \varphi(1)).$ So
$\delta_{f_2}(\varphi(x), \varphi(y)\leq {\rm const}\cdot
\delta_{f_1}(x,y).$ Thus the map $\varphi$ extends to a Lipschitz
map ${\overline \G}_{1, f_1}\to{\overline \G}_{2, f_2}.$\bx

\bigskip

\noindent{\bf Remark.} For the Floyd scaling function
 $\displaystyle f_1={1\over \P_1}$ and for the distortion
function $\al=\P_2$ where $\P_i$ are  polynomials of degree $>1$
we put $\displaystyle f_2=f_1\circ\al={1\over \P_1\circ\P_2}$.
Then the condition (13) is satisfied.

\bigskip

\bigskip

\noindent For a distortion function $\al$ we introduce similarly
to \ref{hull} an $\alpha$-hull $\s H_\alpha M$ of a set $M:$

 $$\displaystyle{\s H_\alpha M}=M{\cup}\{\gamma(I)\ \vert\ \gamma:I\to G\
{\rm is\ a\ \alpha-geodesic},\ I\subset \Z,\ {\rm and}\
\gamma(\partial I) \subset M\}.$$

\noindent We have the following generalization of the Main Lemma
\bigskip

\begin{lem}\label{mlalpha} (Main Lemma for  $\al$-hulls)
For every  polynomial function $\alpha:\N\to \to \R_{\geq 0}$
$$\forall M\subset \ti
T\ \ :\  T \cap \overline M=T \cap \overline{\s H_\al
M}.\hfill\eqno(14)$$

\end{lem}

\proof  Fix  a    Floyd function  $f$ of polynomial type
$f=1/{\cal P}$ such that $\P$ is a polynomial of degree $> 1$ and
the pair $(f,\al)$ is admissible. By Corollary \ref{polynmap} for
the function $f$ there exists the Floyd map $F:\OG\to\ti T$. It
gives the shortcut metric ${\overline \d}_{1,f}$ on $\ti T$
satisfying ${\d}_{1,f}\geq {\overline \d}_{1,f}$ (see (8) and
(8')).

  Using the Generalized Karlsson
Lemma (instead of Karlsson Lemma \ref{agraph}) in the argument of
the Main Lemma we obtain that the Floyd length $\delta_{1,f}$ of
every $\alpha$-geodesic outside of the ball of a radius $r$ is
less that $\ve={\overline \d}_{1,f}(\overline M, a)/2>0$. By the
same argument  it follows that every $\al$-geodesic with endpoints
in $M$ belongs to the union of $\overline M$ and a ball of finite
radius centered at $1\in G.$ The Lemma follows. \bx
\bigskip

\noindent {\bf Definition.} {\it A bi-infinite $\alpha$-geodesic
$\gamma:\Z\to G$ is called {\it  $\alpha$-horocycle} at a point
$p\in T$ if
 $\ga(+\infty)=\gamma(-\infty)=p.$}

 \bigskip

\noindent Using the Generalized Karlsson Lemma  and Lemma
\ref{mlalpha} (instead of  \ref{agraph} and Main Lemma) in the
proofs of Lemmas \ref{hulfin} and \ref{nohorcon} we obtain the
following generalizations for $\al$-geodesics.
\medskip

\begin{lem}
\label{alnohorcon} For every polynomial distortion function $\al$
there is no $\al$-horocycle at a conical point.\bx
\end{lem}

\begin{lem}
\label{alhulfin}

For every polynomial distortion function $\al$ the set $(G \cap \s
H_\al p)/H$ is finite.\bx
\end{lem}

\bigskip

 \noindent {\it Proof of
Theorem C.}   Let $\varphi:H\to G$ be a $\al$-isometric map.
Following the above Remark we fix a Floyd scaling functions $f_1$
and $f_2$  of polynomial type such that $\displaystyle f_1={1\over
\P_1}$ and $\displaystyle f_2={1\over \P_1\circ\al}$ (${\rm deg}
\P_1
> 1$).  By   keeping the  notations of the
proof of Theorem B we put

$$X=\partial_{f_2}G,\ Y=\partial_{f_1}H,\  \ti X =X\sqcup G,\ \ti Y=Y\sqcup H.$$

\noindent  For the Floyd function $f_2$   by  Corollary
\ref{polynmap} there exists the Floyd map $F:\ti X\to   \ti T.$

As in the proof of Theorem B we define  the quotients $\ti S$ and
$S$ of the spaces $\ti Y$ and $Y$ respectively. The action
$H\act\ti S$ is 3-discontinuous. To show that it is 2-cocompact we
modify the proof of  Proposition \ref{cocomppairs} as follows. For
a geodesic $\ga\subset H$ we obtain the $\al$-geodesic
$\varphi(\ga)\subset G.$

The series $\displaystyle \Sigma=\sum_k\al_{2k+1}f_{2,k}$
converges. Indeed  $\al$ is a polynomial function so
$$\displaystyle{\al(2k+1)\over \P_1(\al(k))}\leq {\rm const}\cdot
{\al(k)\over\P_1(\al(k))}.$$

\noindent Since ${\rm deg} \P_1 > 1$ the series $\displaystyle
\sum_k{\al(k)\over \P_1(\al(k))}$ converges and so is $\Sigma.$
Thus the Generalized Karlsson Lemma \ref{genkarl} applies to the
space $\ti X$. So we conclude as before that $\varphi(\ga)$ is
situated within a distance $r=r(\al, \delta)$ from the basepoint
$g.$ So $d(g, g_0)\leq r$ where $g_0=\varphi(h_0)$ for some
$h_0\in H.$ Applying now Lemma \ref{alhomfl} (instead of
\ref{homfl}) we conclude that that the action $H\act\Th^2S$ admits
a compact fundamental set.

To show that $\varphi_*(p)\in T$ is parabolic for a parabolic
point $p\in S$  we proceed similarly. Since every parabolic
subgroup is undistorted in $H$ there exists a $c$-quasigeodesic
horocycle $\ga\subset H$ at the limit point $p$. The composition
$\varphi\circ\ga$ is an $\al$-geodesic horocycle in $\ti T.$  It
follows from Lemma \ref{alnohorcon}  that its unique limit point
$\varphi_*(p)\in T$ is parabolic for the action $G\act T.$

Similarly applying  Lemma \ref{alhulfin}    we obtain that the set
$G\cap H_{\al}q$ is ${\rm Stab}_Gq$-finite. It implies that the
image of every parabolic subgroup $Q$ of $H$ is situated in a
uniformly bounded distance from a right coset of
$\s{Stab}_G(\varphi_* p)$ in $G.$ Theorem C follows.\bx

\section{Appendix: a short proof that 3-cocompactness of an action
implies  word-hyperbolicity of the group.}

As an application of our method we give a short proof of the
following theorem of B.~Bowditch:

\bigskip

\noindent{\bf Theorem} \cite{Bo3}. {\it Let $G$ be a finitely
generated group acting $3$-discontinuously and $3$-cocompactly on
a compactum $T$ without isolated points. Then $G$ is
word-hyperbolic.}

\medskip

\noindent The following lemma requires some additional information
from [Ge1].

\medskip

\begin{lem}
\label{nonConicPair} Let a group $G$ act $3$-discontinuously on a
compactum $T$. Let ${p,q}$ be distinct non-conical points in $T$
and let $K$ be a compact subset of $T^2{\setminus}\De^2T$. Then
the set\hfil\penalty-10000 $S{=}\{g\in G:(gp,gq)\in K\}$ is
finite.\end{lem}

 \proof Assume that $S$ is infinite. The compact
$K$ can be covered by finitely many closed subproducts of the form
$A{\times}B$ with $A{\cap}B{=}\varnothing$. So we can assume that
$K{=}A{\times}B$ where $A,B$ are closed disjoint sets. The set
$\Lambda_{\s{rep}}S$ of the repellers of the limit crosses for $S$
(see
 [Ge1, subsection 18]) is nonempty. It is contained in
$\{{p,q}\}$ since otherwise, for some $g\in S$, the pair
$\{gp,gq\}$ becomes arbitrarily small.

So $S$ contains an infinite subset $S_1$ with
$\Lambda_{\s{rep}}S_1$ being a single point. Without loss of
generality we can assume that $S{=}S_1$ and
$\Lambda_{\s{rep}}S{=}\{ p\}$. The set $\Lambda_{\s{attr}}S$ of
the attractors of the limit crosses is contained in $B$. Let $B_1$
be a closed neighborhood of $B$ disjoint from $A$. Thus, for $
a\in T{\setminus}\{p\}$, the set $\{g\in S:ga\notin  B_1\}$ is
finite since it possesses no limit crosses. Hence
$\{\{gp,ga\}:g\in S\}$ is contained in a compact subset of
$\Th^2T$. So $ p$ is conical by  of [Ge1, Definition 3]. \bx

\bigskip

\noindent {\bf Remark.} With an additional assumption that $T$ is
metrisable the  Lemma can be also proved  using Gehring-Martin's
definition of the convergence property.

\medskip

\noindent \sc Corollary\sl. If $G$ acts $3$-discontinuously and
$3$-cocompactly on a compactum $T$ without isolated points then
every point of $T$ is conical.

\medskip

\noindent  \it Proof\rm.
 Clearly, the $3$-cocompactness implies the $2$-cocompactness, hence every nonconical
  point is either isolated or parabolic \cite{Ge1}.
By the assumption the discontinuity domain is empty. Assume that
parabolic points do exist. Hence there exist at least two
parabolic points since otherwise we must have a discontinuity
domain.

Let ${p,q}$ be distinct parabolic points and let $L$ be a compact
fundamental set for the action of $G$ on $\Th^3T$. We can assume
that $L$ has the form
 $\bigcup\limits_{i=1}^nA_i{\times}B_i{\times}C_i$ where $A_i, B_i, C_i$
  are closed disjoint subsets of $T.$
 For every $a\in T{\setminus}\{{p,q}\}$ there exist $g_{ a}\in G$ such that
$g_{ a}(p,q,a)\in L$. By Lemma \ref{nonConicPair} the set
$\{g_{a}: a\in T{\setminus}\{{p,q}\}\}$ is finite. If $g_{
a}(p,q,a)\in A_{i(a)}{\times}B_{i(a)}{\times}C_{i(a)}$ then
$T{\setminus}\{{p,q}\}$ is a union of finitely many closed sets
$g_a^{-1}C_{i(a)}$. Thus $\{p,q\}$ is open and $p$ and $q$ are
isolated contradicting the assumption.\bx

\bigskip

\noindent\it Proof of the Theorem\rm. Assume that $G$ is not
hyperbolic. There exists a sequence of
 geodesic triangles  with the sides $\{l_n, m_n, k_n\}$  so that
$d(x_n, m_n\cup k_n)\to\infty$ for $x_n\in l_n.$
Using the $G$-action we can make $x_n$ equal to $1$ for all $n$. By
Karlsson Lemma the Floyd length of $m_n\cup k_n$ tends to zero and so $\d_1(y_n, z_n)\to
0$ where $\partial l_n=\{y_n, z_n\}.$ Since all $l_n$ pass through
the same point $1$ we can choose  a subsequence converging to
a geodesic horocycle $l.$ By Lemma \ref{nohorcon} the target of $l$
is not conical contradicting the above Corollary.  \bx

\medskip

\noindent {\bf Remark.} In the paper \cite{Bo3} the above Theorem
is proved without the assumption that  $G$ is finitely generated.
It is shown in [GePo2, Corollary 3.38] that a group $G$ acting
$3$-discontinuously and $2$-cocompactly on a compactum $T$ is
relatively finitely generated with respect to the maximal
parabolic subgroups. In particular if $G$ acts 3-cocompactly on
$T$ by the above Corollary there are no parabolic subgroups and so
$G$ is finitely generated.

\noindent Victor Gerasimov : Department of Mathematics, Federal
University of Minas Gerais,

\noindent Av. Antônio Carlos, 6627 - Caixa Postal 702 - CEP
30161-970 Belo Horizonte, Brazil;

\noindent email: victor@mat.ufmg.br

\bigskip

\noindent Leonid Potyagailo: UFR de Math\'ematiques, Universit\'e
de Lille 1,

\noindent 59655 Villeneuve d'Ascq cedex, France;

\noindent email: Leonid.Potyagailo@univ-lille1.fr


\bigskip



\begin{thebibliography}{Main25}
\addcontentsline{toc}{chapter}{Bibli ography}

\bibitem[BM]{BM}

 A. Beardon, B. Maskit, {\it Limit sets of Kleinian groups and
finite sided fundamental polyhedra}, Acta Math. \bf 132, \rm 1974,
1--12.

\bibitem[BDM]{BDM}
J.~Behrstock, C.~Drutu, L.~Mosher, {\it Thick metric spaces,
relative hyperbolicity, and quasi-isometric rigidity}, Math.
Annalen, {\bf 344}, 2009, 543-595.



\bibitem[Bo1]{Bo1} B. H. Bowditch, {\it Relatively hyperbolic
groups,}  1997, Preprint.

\bibitem[Bo2]{Bo2} B. H. Bowditch, {\it Convergence groups and configuration
spaces},  in ``Group theory down under'' (ed.\ J.Cossey,
C.F.Miller, W.D.Neumann, M.Shapiro), de Gruyter (1999) 23--54.

\bibitem[Bo3]{Bo3}
 B. H. Bowditch, {\it A topological characterisation of hyperbolic
groups}, J. Amer. Math. Soc. 11, 1998, no. 3, 643--667.

\bibitem[BH]{BH}
M.~Bridson, A.~Haefliger, {\it Metric spaces of non-positive
curvature}, Grundlehren der Mathematischen Wissenschaften, 319.
Springer-Verlag, Berlin, 1999.

\bibitem[Bourb]{Bourb}
 N.~Bourbaki, {Topologie G\'en\'erale}
Hermann, Paris, 1965.

\bibitem[BBI]{BBI}
D.~Burago, Y.~Burago, S.~Ivanov, {A Course in metric geometry}
Graduate Studies in Math, AMS, Vol 13, 2001.




\bibitem[Dr]{Dr}
 C.~Drutu, {\it Relatively hyperbolic groups: geometry and quasi-isometric
invariance},  Comment. Math. Helv., {\bf 84} (3), 2009, 503-546.




\bibitem[Fa]{Fa} B.~Farb, {\it Relatively hyperbolic groups,}
GAFA, 8(5), 1998, 810-840.

\bibitem[F]{F}
W. J. Floyd, {\it Group completions and limit sets of Kleinian
groups}, Inventiones Math. \bf 57, \rm 1980, 205--218.


\bibitem[Fr]{Fr}
 E. M. Freden, {\it Properties of convergence groups and spaces}, Conformal
Geometry and Dynamics, \bf 1, \rm 1997, 13--23.


\bibitem[Ge1]{Ge1}
V.~Gerasimov, {\it Expansive convergence groups are relatively
hyperbolic},  GAFA, {\bf 19}, 2009, 137--169.

\bibitem[Ge2]{Ge2}
V.~Gerasimov, {\it Floyd maps to the boundaries of relatively
hyperbolic groups}, arXiv:1001.4482 [math.GR], 2010.


\bibitem[GePo1]{GePo1}
V.~Gerasimov, L.~Potyagailo, {\it Quasi-isometric maps and Floyd
boundaries of relatively hyperbolic groups}, arXiv:0908.0705
[math.GR], 2009.


\bibitem[GePo2]{GePo2}
V.~Gerasimov, L.~Potyagailo, {\it Horospherical geometry of
relatively hyperbolic groups}, arXiv:1008.3470[math.GR], 2010.


\bibitem[GePo3]{GePo3} V.~Gerasimov, L.~Potyagailo, {\it
Quasiconvexity in the relatively hyperbolic groups}, In
preparation.

\bibitem[GM]{GM}
F. W. Gehring and G. J. Martin\sl, Discrete quasiconformal groups
I\rm. Proc.\ London Math\ Soc.,\ \bf 55, \rm 1987, 331--358.

\bibitem[Gr]{Gr} M. Gromov, {\it Hyperbolic groups}, in: ``Essays in Group
Theory'' (ed. S.~M.~Gersten) M.S.R.I. Publications No.~8,
Springer-Verlag (1987) 75--263.

\bibitem[Gr1]{Gr1} M. Gromov, {\it Asymptotic invariants of infinite groups},
``Geometric Group Theory II'' LMS Lecture notes \bf 182\rm,
Cambridge University Press (1993)

\bibitem[Hr]{Hr}
G.~Hruska, {\it Relative hyperbolicity and relative quasiconvexity
for countable groups}, arXiv:0801.4596 [math.GR], 2008.


\bibitem[Ka]{Ka}
 A. Karlsson {\it Free subgroups of groups with non-trivial Floyd boundary}, Comm.
Algebra, 31, 2003, 5361--5376.

\bibitem[Ke]{Ke}
J.L. Kelly, {\it General Topology}, Graduate Texts in Mathematics,
N 27, Springer Verlag, New York, 1975.

\bibitem[My]{My}
P.~J. Myrberg {\it Untersuchungen ueber die automorphen Funktionen
beliebiger vieler Variabelen}, Acta Math. \bf 46, \rm 1925.


\bibitem[Os]{Os}
D.~Osin, {\it Relatively hyperbolic groups: intrinsic geometry,
algebraic properties and algorithmic problems}, \rm, Mem. AMS \bf
179 \rm(2006) no. 843 vi+100pp.



\bibitem[Tu1]{Tu1}
P. Tukia, {\it A remark on the paper by Floyd}, Holomorphic
functions and moduli, vol.II(Berkeley CA,1986),165--172, MSRI
Publ.11, Springer New York 1988.

\bibitem[Tu2]{Tu2}
  P. Tukia, {\it Convergence groups and Gromov's metric
hyperbolic spaces} : New Zealand J.\ Math.\ \bf 23, \rm 1994,
157--187.

\bibitem[Tu3]{Tu3}
 P. Tukia, {\it Conical limit points and uniform convergence
groups}, J.\ Reine.\ Angew.\ Math.\ {\bf 501}, 1998, 71--98.

\bibitem[W]{W}
A.~Weil, {\rm Sur les espace à structure uniforme et sur la
topologie générale}, Paris, 1937.

\bibitem[Ya]{Ya}
A. Yaman, {\it A topological characterisation of relatively
hyperbolic groups}, J.\ reine ang.\ Math. \bf 566, \rm 2004,
41--89.




\end{thebibliography}
\end{document}